\documentclass[10pt]{article} 
\usepackage[utf8]{inputenc}
\usepackage{geometry}                      
\usepackage{graphicx}   
\usepackage{xcolor}             
\usepackage{amssymb}
\usepackage{amsmath}
\usepackage{amsthm}

\newtheorem{exmp}{Example}

\newtheorem{fact}{Fact}
\newtheorem{rmrk}{Remark}

\newtheorem{algo}{Algorithm}
\usepackage[]{hyperref}
\hypersetup{
    colorlinks = true,
    linkbordercolor = {white},
    linkcolor=black,
    citecolor=black,
    urlcolor=black
}

\begin{document}
\title{Approximating the Perfect Sampling Grids for Computing the Eigenvalues of Toeplitz-like Matrices Using the Spectral Symbol}
\author{Sven-Erik Ekström\\ {\small\texttt{see@2pi.se}}}
\date{Athens University of Economics and Business\\[2ex]%
    \today
}
\maketitle
\begin{abstract}
In a series of papers the author and others have studied an asymptotic expansion of the errors of the eigenvalue approximation, using the spectral symbol, in connection with Toeplitz (and Toeplitz-like) matrices, that is, $E_{j,n}$ in  $\lambda_j(A_n)=f(\theta_{j,n})+E_{j,n}$, $A_n=T_n(f)$, $f$ real-valued cosine polynomial.
In this paper we instead  study an asymptotic expansion of the errors of the equispaced sampling grids $\theta_{j,n}$, compared to the exact grids $\xi_{j,n}$ (where $\lambda_j(A_n)=f(\xi_{j,n})$), that is, $E_{j,n}$ in $\xi_{j,n}=\theta_{j,n}+E_{j,n}$.
We present an algorithm to approximate the expansion. Finally we show numerically that this type of expansion works for various kind of Toeplitz-like matrices (Toeplitz, preconditioned Toeplitz, low-rank corrections of them). We  critically discuss several specific examples and we demonstrate the superior numerical behavior of the present approach with respect to the previous ones.
\end{abstract}

\section{Introduction}
\label{sec:introduction}
For a given banded Hermitian Toeplitz matrix, $T_n(f)\in \mathbb{R}^{n\times n}$,
\begin{align}
T_n(f)=\left[
\begin{array}{ccccccccc}
\hat{f}_0&\cdots&\hat{f}_m\\
\vdots&\ddots&&\ddots\\
\hat{f}_m&&\ddots&&\ddots\\
&\ddots&&\ddots&&\hat{f}_m\\
&&\ddots&&\ddots&\vdots\\
&&&\hat{f}_m&\cdots&\hat{f}_0
\end{array}
\right]\nonumber
\end{align}
we can easily associate a real-valued spectral symbol $f$, independent of $n$, namely
\begin{align}
f(\theta)=\sum_{k=-n}^{n}\hat{f}_ke^{\mathbf{i}k\theta}=\sum_{k=-m}^{m}\hat{f}_ke^{\mathbf{i}k\theta}=\hat{f}_0+2\sum_{k=0}^{m}\hat{f}_k\cos(\theta),\nonumber
\end{align}
and note that $\hat{f}_k$ are the Fourier coefficients of $f(\theta)$, that is,
\begin{align}
\hat{f}_k=\frac{1}{2\pi}\int_{-\pi}^{\pi}f(\theta)e^{-\mathbf{i}k\theta}\mathrm{d}\theta.\nonumber
\end{align}
We say that the symbol $f$ generates the Toeplitz matrix of order $n$ when  $T_n(f)$ is defined as $[\hat{f}_{i-j}]_{i,j=1}^n$.
The eigenvalues of $T_n(f)$, denoted $\lambda_j(T_n(f))$, can be approximated by sampling the symbol $f$ with an equispaced grid $\theta_{j,n}$, that is,
\begin{align}
\lambda_j(T_n(f))=f(\theta_{j,n})+E_{j,n}^{\lambda,\theta}
\label{eq:introduction:lambdathetae0}
\end{align}
where we typically use the following standard equispaced grid
\begin{align}
\theta_{j,n}=\frac{j\pi}{n+1}=j\pi h,\quad j=1,\ldots,n,\quad h=\frac{1}{n+1},
\label{eq:introduction:thetagrid}
\end{align}
and the error, for grid $\theta_{j,n}$, is $E_{j,n}^{\lambda,\theta}=\mathcal{O}(h)$; see for example~\cite{bottcher991,garoni171}.

\begin{figure}[!ht]
\centering
\includegraphics[width=0.48\textwidth]{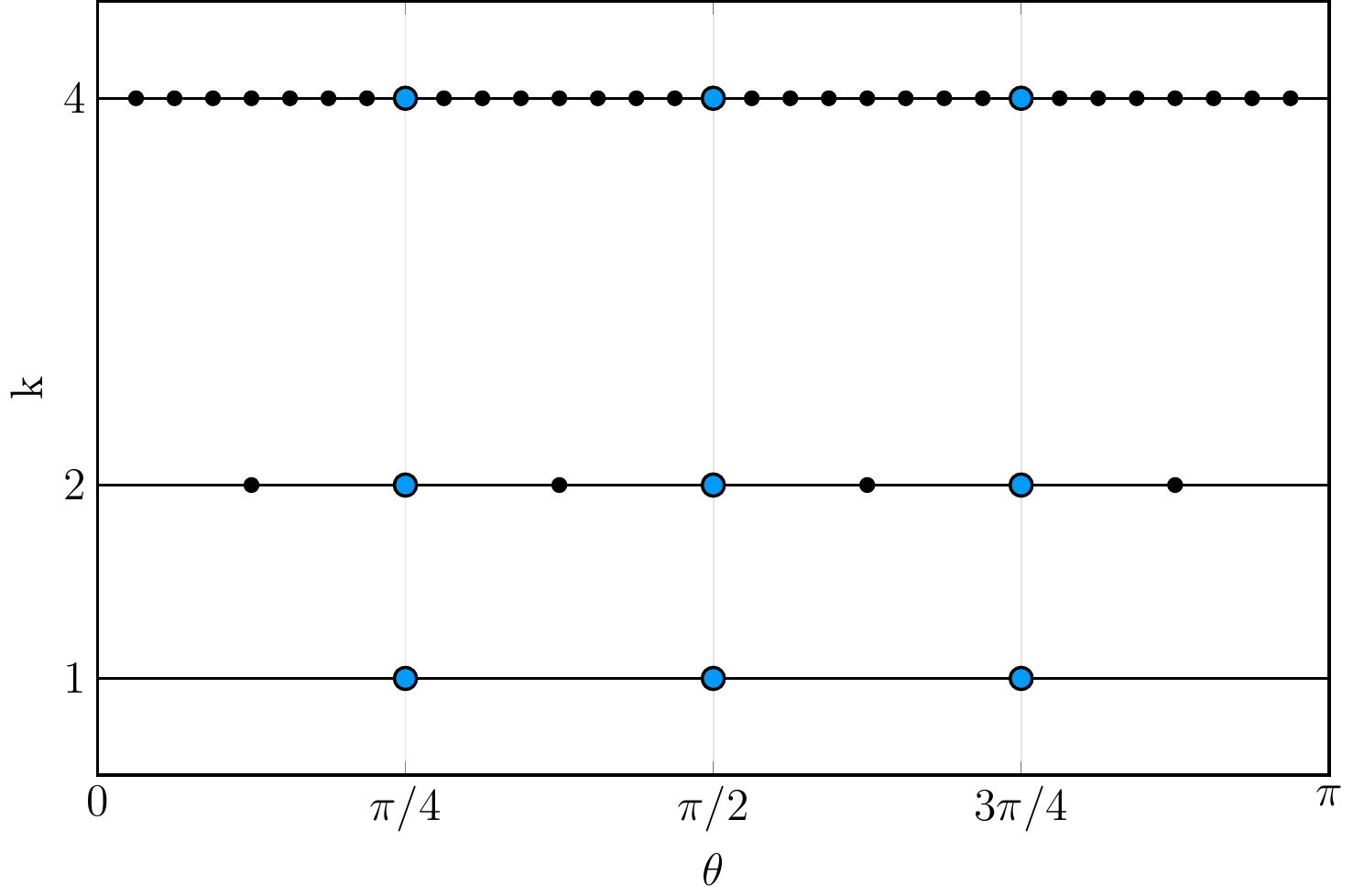}
\caption{Three grids $\theta_{j,n_k}$ with $k=1,2,4$, $n_k=2^{k-1}(n_1+1)-1$, and $n_1=3$. The subgrids $\theta_{j_k,n_k}$, which is the same for any $k$, is indicated by blue circles ($j_k=2^{k-1}j_1$, $j_1=\{1,\ldots, n_1\}$).}
\label{fig:introduction:grid}
\end{figure}

In a series of papers~\cite{bogoya151,bogoya171,bottcher101} an asymptotic expansion of
the error $E_{j,n}^\theta$ in \eqref{eq:introduction:lambdathetae0},
\begin{align}
E_{j,n}^{\lambda,\theta}=\sum_{k=1}^\alpha c_k(\theta)h^k+E_{j,n,\alpha}^{\lambda,\theta}
\label{eq:introduction:errorexpansiontheta}
\end{align}
was studied. In~\cite{ekstrom171} an algorithm to approximate the functions $c_k(\theta)$, in the grid points $\theta_{j,n_1}$ for a chosen $n_1$, was proposed: such an approximation is denoted by $\tilde{c}_k(\theta_{j,n_1})$, for $k=1,\ldots,\alpha$ and $\alpha$ chosen properly.
We thus have
\begin{align}
E_{j,n_1}^{\lambda,\theta}&=\sum_{k=1}^\alpha c_k(\theta_{j,n_1})h^k+E_{j,n_1,\alpha}^{\lambda,\theta}\nonumber\\
&=\sum_{k=1}^\alpha \tilde{c}_k(\theta_{j,n_1})h_1^k+\tilde{E}_{j,n_1,\alpha}^{\lambda,\theta}.\nonumber
\end{align}
Ignoring the error term $\tilde{E}_{j,n,\alpha}^{\lambda,\theta}$, we obtain the following expression to approximate the eigenvalues of a generated matrix $T_{n_m}(f)$ with  indices $j_m=2^{m-1}j_1$ (where $j_1=\{1,\ldots,n\}$), that is
\begin{align}
\lambda_{j_m}(T_{n_m}(f))\approx\tilde{\lambda}_{j_m}^{(\alpha)}(T_{n_m}(f))= f(\theta_{j,n_1})+\sum_{k=1}^\alpha \tilde{c}_k(\theta_{j,n_1})h_m^k,\nonumber
\end{align}
where $n_m=2^{m-1}(n_1+1)-1$. In Figure~\ref{fig:introduction:grid} we show three grids, given by $n_k=2^{k-1}(n_1+1)-1$, namely $n_1=3, n_2=7$, and $n_4=31$. Blue circles indicate the sets of incices $j_k$ such that $\theta_{j_k,n_k}$ always are the same for any $k$.

\noindent If we use the interpolation--extrapolation scheme introduced in \cite{ekstrom183}, then we can approximate $\tilde{c}_k(\theta_{j,n})$ for any $n$ and thus we can compute numerically the full spectrum of $T_n(f)$ by
\begin{align}
\lambda_{j}(T_{n}(f))\approx\tilde{\lambda}_{j}^{(\theta,\alpha)}(T_{n}(f))= f(\theta_{j,n})+\sum_{k=1}^\alpha \tilde{c}_k(\theta_{j,n})h^k.
\label{eq:introduction:lambdaexpansiontheta}
\end{align}
Note that the procedure using the asymptotic expansion \eqref{eq:introduction:lambdaexpansiontheta} also works for more general types of Toeplitz-like matrices, for example for preconditioned Toeplitz matrices~\cite{ahmad171}, and block-Toeplitz matrices generated by matrix-valued symbols~\cite{ekstrom181}, and even in a differential setting when considering the isogeometric approximation of the Laplacian eigenvalue problem in any dimension \cite{ekstrom185}. The latter example is quite important since the involved matrices can be written as low-rank corrections of either of Toeplitz or of preconditioned Toeplitz matrices. We call this type of methods \textit{matrix-less} since they do not need to construct a matrix, of arbitrary order $n$, to approximate its eigenvalues (after the initial step when approximating $c_k(\theta_{j,n_1})$ from $\alpha$ small matrices, data that can be reused for a matrix of any order).

By the notation of the theory of Generalized Locally Toeplitz (GLT) sequences, see \cite{garoni171} and references therein, we say that $f$ describes the eigenvalue distribution of a matrix $A_n$, when
\begin{align}
\{A_n\}_n\sim_{\textsc{glt},\sigma,\lambda}f.\nonumber
\end{align}
Here we do not report the formal definitions for which we refer to  \cite{garoni171}, but instead we remind that
$\{A_n\}_n\sim_{\sigma,\lambda}f$ has the informal meaning that the bulk of singular values, eigenvalues behave, for $n$ large enough and up to infinitesimal errors in the parameter $n$, as an equispaced sampling of the functions $f$, $|f|$, respectively over the common definition domain.

In this paper we consider the following types of matrices $A_n$,
\begin{align}
A_n&=T_n(f),\nonumber\\
A_n&=T_n(f)+R_n,\nonumber\\
A_n&=T_n^{-1}(b)T_n(a),\quad f=a/b,\nonumber
\end{align}
where $f$ is assumed to be real-valued, continuous and monotone, $R_n$ is a low-rank matrix, and $a,b$ are real-valued cosine polynomials, with $b$ non-negative and not identically zero. We thus have the following more general notation
\begin{align}
\tilde{\lambda}_{j}^{(\theta,\alpha)}(A_n)&= f(\theta_{j,n})+\sum_{k=1}^\alpha \tilde{c}_k(\theta_{j,n})h^k,\label{eq:introduction:lambdaaproxexpansiontheta}\\
\tilde{E}_{j,n,\alpha}^{\lambda,\theta}&=\lambda_{j}(A_n)-\tilde{\lambda}_{j}^{(\alpha)}(A_n).
\label{eq:introduction:errorapproxexpansiontheta}
\end{align}

\section{Main Results}
\label{sec:main}
In this paper we now introduce the ``perfect'' grid $\xi_{j,n}$, associated with the symbol $f$ and matrix $A_n$, such that
 \begin{align}
\lambda_j(A_n)=f(\xi_{j,n}),
\label{eq:main:lambdaxi}
\end{align}
that is, sampling the symbol $f(\theta)$ with $\xi_{j,n}$ gives the exact eigenvalues of $A_n$. Such a grid has to exist thanks to the continuity of $f$ and thanks to the distributional/localization results (see \cite{dibenedetto931} and references therein). We thus have
\begin{align}
\xi_{j,n}=\theta_{j,n}+E_{j,n}^\xi,
\label{eq:main:lambdaxie0}
\end{align}
where $E_{j,n}^\xi$ is error in the sampling grid, using the standard grid~\eqref{eq:introduction:thetagrid} instead of the perfect grid $\xi_{j,n}$.
We now introduce the following fact, and dedicate the rest of this paper to show applications of this assumption to a number of different Toeplitz-like matrices and their symbols, describing the eigenvalue distribution. Fact~\ref{fact:main:1} follows from \eqref{eq:introduction:errorexpansiontheta}  and from the monotonicity of $f$ via a proper Taylor expansion.

\begin{fact}
\label{fact:main:1}
There exists an expansion 
\begin{align}
\xi_{j,n}=\theta_{j,n}+\sum_{k=1}^\alpha d_k(\theta_{j,n})h^k+E_{j,n,\alpha}^\xi,
\label{eq:main:xiexpansion}
\end{align}
where $h=1/(n+1)$,
such that we can asymptotically describe a perfect grid $\xi_{j,n}$ such that \eqref{eq:main:lambdaxi} is true. Of course the latter holds for the same types of symbols for which the expansion~\eqref{eq:introduction:errorexpansiontheta} has been shown to work for; see for example~\cite{ahmad171,ekstrom181,ekstrom171}.
\end{fact} 
We now first present two illustrative examples, and then propose Algorithm~\ref{algo:main:1}, a matrix-less method for approximating the functions $d_k(\theta)$ in \eqref{eq:main:xiexpansion}. Finally Algorithm~\ref{algo:main:2} describes how to use the information form Algorithm~\ref{algo:main:1} to approximate the eigenvalues $\lambda_j(A_n)$ for any $n$.
Henceforth we use the following notation
\begin{align}
\tilde{\xi}_{j,n}^{(\alpha)}&=\theta_{j,n}+\sum_{k=1}^\alpha \tilde{d}_k(\theta_{j,n})h^k,\label{eq:main:xiapproxexpansion}\\
\tilde{E}_{j,n,\alpha}^\xi&=\xi_{j,n}-\tilde{\xi}_{j,n}^{(\alpha)},\label{eq:main:errorxiapproxxiexpansion}\\
\tilde{\lambda}_{j}^{(\xi,\alpha)}(A_n)&=f(\xi_{j,n}^{(\alpha)}),\label{eq:main:lambdaapproxxiexpansion}\\
\tilde{E}_{j,n,\alpha}^{\lambda,\xi}&=\lambda_j(A_n)-\tilde{\lambda}_{j,n}^{(\xi,\alpha)}(A_n),\label{eq:main:errorlambdaapproxxiexpansion}
\end{align}
to distinguish from the expansion and errors defined in \eqref{eq:introduction:lambdaaproxexpansiontheta} and \eqref{eq:introduction:errorapproxexpansiontheta}.
\begin{exmp}
\label{exmp:main:errorlambdavsxi} 
In this example we study the symbol $f(\theta)=(2-2\cos(\theta))^2=6-8\cos(\theta)+2\cos(2\theta)$, associated with the second order finite difference discretization of the bi-Laplacian. It is well known that an equispaced grid does not give the exact eigenvalues of the generated matrices $T_n(f)$ for finite matrices, for example, see discussions in \cite{barrera181,ekstrom184}.

In Figure~\ref{fig:main:symbol} we show the symbol $f$ (dashed pink line), and present data for the generated Toeplitz matrix of order $n=3$, $T_3(f)$, that is, the three eigenvalues $\lambda_j(T_3(f))=f(\xi_{j,3})$ (red circles) and the three samplings $f(\theta_{j,3})$ (blue circles) where the sampling grid $\theta_{j,3}$ is defined in~\eqref{eq:introduction:thetagrid}.

\begin{figure}[!ht]
\centering
\includegraphics[width=0.48\textwidth]{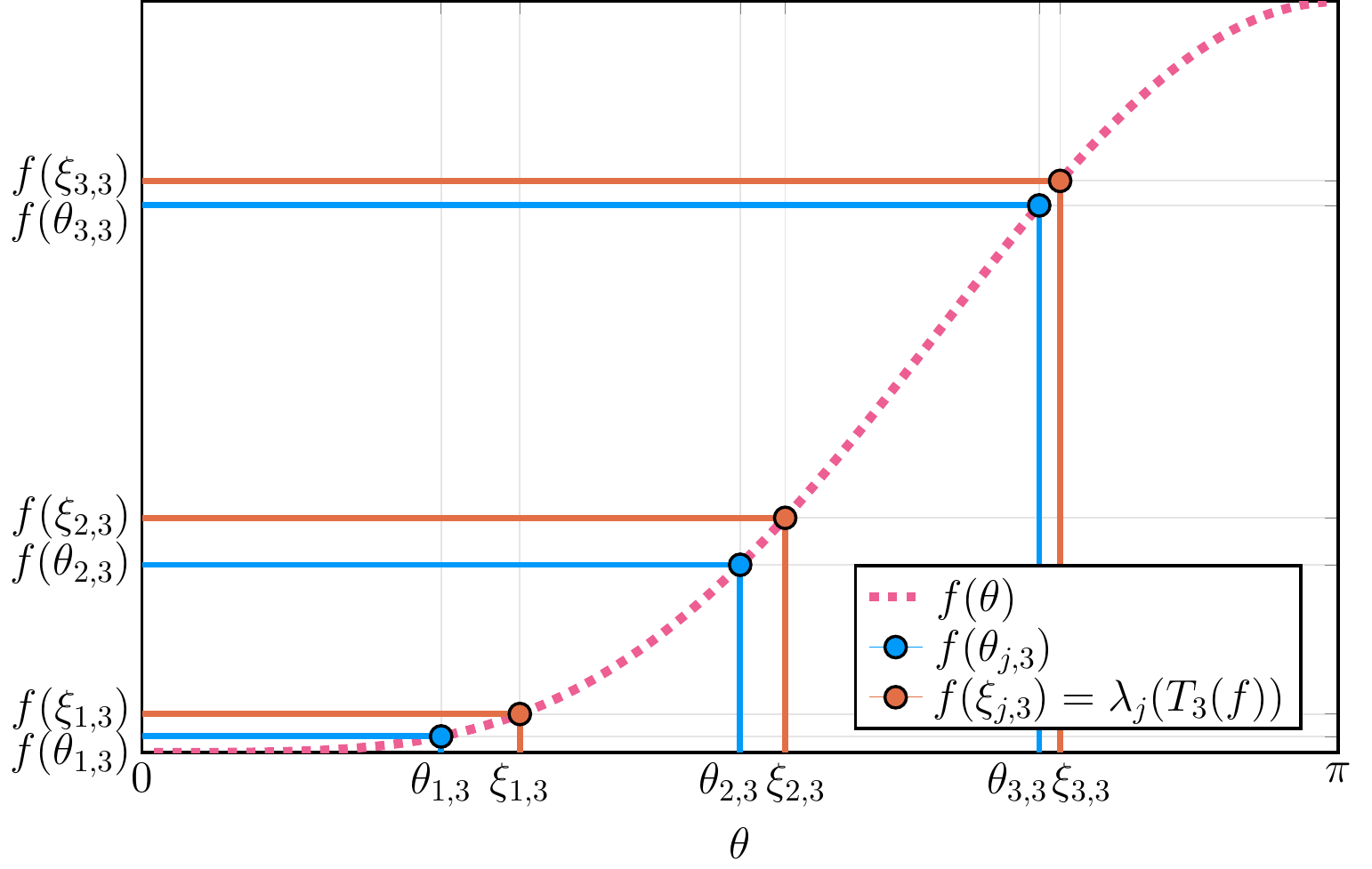}
\caption{Example~\ref{exmp:main:errorlambdavsxi}: Symbol $f(\theta)=(2-2\cos(\theta))^2$ (dashed pink line) and the two expansions approximate the errors $E_{j,n}^{\lambda,\theta}=\lambda_j(T_n(f))-f(\theta_{j,n})=f(\xi_{j,n})-f(\theta_{j,n})$ and $E_{j,n}^\xi=\xi_{j,n}-\theta_{j,n}$.}
\label{fig:main:symbol}
\end{figure}

Hence, in previous papers, we approximate $E_{j,n}^{\lambda,\theta}$ in \eqref{eq:introduction:lambdathetae0}, by computing the $\tilde{\lambda}_j^{(\theta,\alpha)}(A_n)-f(\theta_{j,n})$ in \eqref{eq:introduction:lambdaaproxexpansiontheta}. We thus approximate the error of the \textit{eigenvalue approximations} when only using the symbol approximation $f(\theta_{j,n})$. The sought error $E_{j,3}^{\lambda,\theta}=\lambda_j(T_3(f))-f(\theta_{j,3})$ is the distance between the red and blue line on the y-axis.

In this paper we instead study the error $E_{j,n}^\xi$ in  \eqref{eq:main:lambdaxie0}, by computing the $\tilde{\xi}_{j,n}^{(\alpha)}-\theta_{j,n}$ in \eqref{eq:main:xiapproxexpansion}. Hence, we approximate the error of the \textit{grid}, when using the standard grid \eqref{eq:introduction:thetagrid}. The sought error $E_{j,3}^\xi=\xi_{j,3}-\theta_{j,3}$ is the distance between the red and blue line on the x-axis.
Subsequently we can use the approximated $\tilde{\xi}_{j,n}^{(\alpha)}$, in \eqref{eq:main:xiapproxexpansion}, to approximate the eigenvalues $\lambda_j(T_n(f))$ by using \eqref{eq:main:lambdaapproxxiexpansion}.
\end{exmp}
 
In the following Example~\ref{exmp:main:errorxi} we show the first numerical evidence supporting Fact~\ref{fact:main:1}.

\begin{exmp} 
\label{exmp:main:errorxi}
We return to the symbol $f(\theta)=(2-2\cos(\theta))^2$ of Example~\ref{exmp:main:errorlambdavsxi}. In Figure~\ref{fig:main:errorxi} we present in the left panel the errors $E_{j,n}^\xi=E_{j,n,0}^\xi$ in \eqref{eq:main:xiexpansion}, for a number of different $n$, that is,
\begin{align}
\xi_{j,n}&=\cos^{-1}\left(\frac{2-\sqrt{\lambda_j(T_n(f))}}{2}\right),\nonumber\\
E_{j,n}^\xi&=\xi_{j,n}-\theta_{j,n}.\nonumber
\end{align}
The sizes of the matrices, used to compute $E_{j,n_k}^\xi$ are $n_k=2^{k-1}(n_1+1)-1$, where $n_1=100$, and $k=1,\ldots,4$. In the right panel we have scaled all the errors by $h_k=1/(n_k+1)$. The overlap is very good for the curves for different $n_k$, indicating the general shape of $d_1(\theta)$ in \eqref{eq:main:xiexpansion}, later observed in Figure~\ref{fig:numerical:bilaplaceexpansion}.

\begin{figure}[!ht]
\centering
\includegraphics[width=0.49\textwidth]{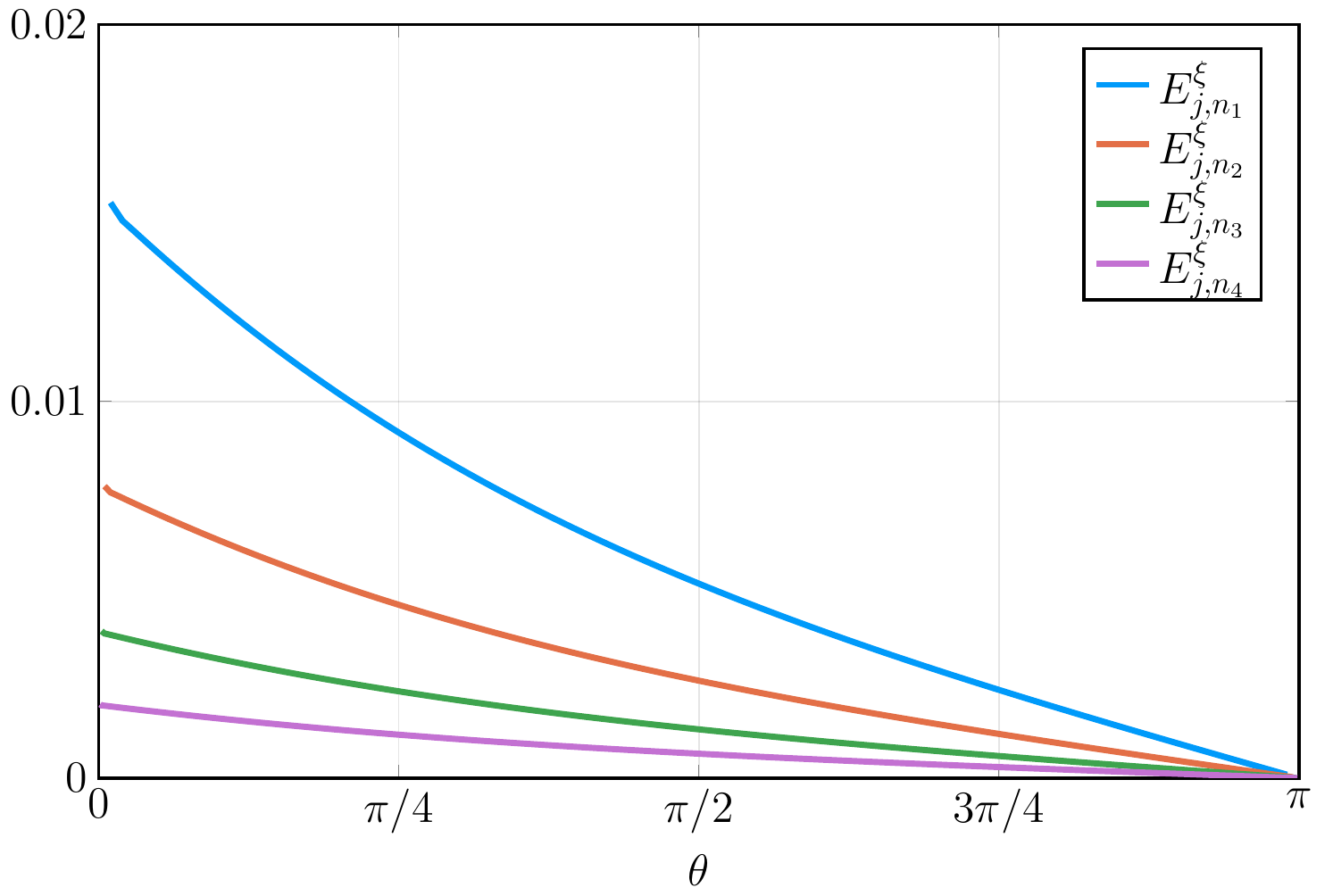}
\includegraphics[width=0.47\textwidth]{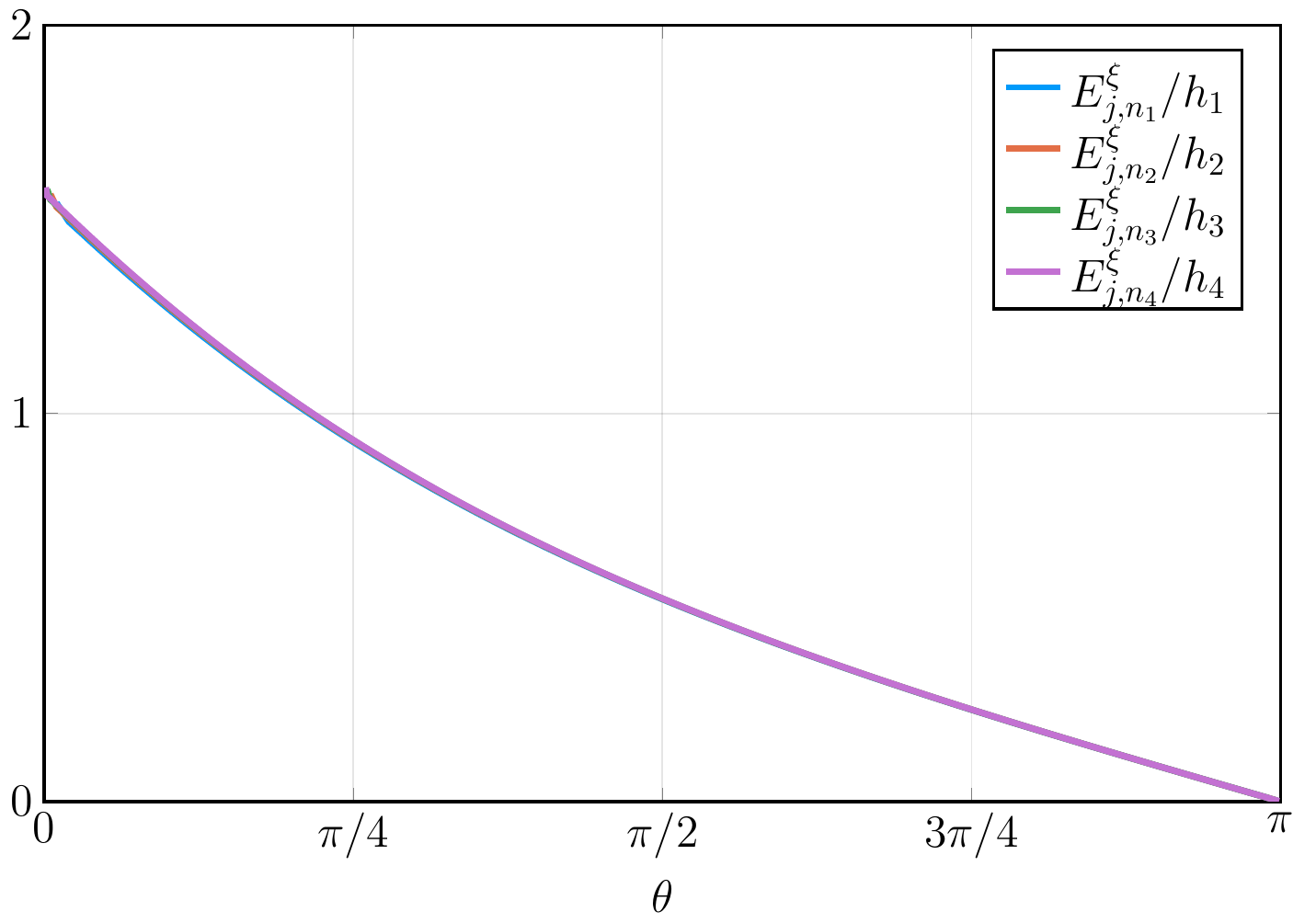}
\caption{Example \ref{exmp:main:errorxi}: Errors $(E_{j,n}^\xi=\xi_{j,n}-\theta_{j,n})$ and scaled errors $(E_{j,n}^\xi/h)$. Left: Errors $E_{j,n_k}^\xi$, for $n_k=2^{k-1}(n_1+1)-1$, $n_1=100$, and $k=1,\ldots,4$. Right: Scaled errors $E_{j,n_k}^\xi/h_k$ where $h_k=1/(n_k+1)$.}
\label{fig:main:errorxi}
\end{figure}
\begin{rmrk}
For more complicated symbols $f$, the roots of $f(\theta)-\lambda_j(A_n)$ are $\xi_{j,n}$, and can be computed numerically. 
\end{rmrk}
\end{exmp}
\noindent We propose in Algorithm~\ref{algo:main:1} a procedure for approximating $d_k(\theta)$ with $\tilde{d}_k(\theta_{j,n_1})$, for $k=1,\ldots,\alpha$ on the grid points $\theta_{j,n_1}$, defined in \eqref{eq:introduction:thetagrid}. In Algorithm~\ref{algo:main:2} we then use the interpolation--extrapolation scheme described in~\cite{ekstrom183} to approximate $d_k(\theta_{j,n})$ for any $n$, and subsequently use this data to approximate the eigenvalues $\lambda_j{A_n}$. We here note that a similar approach was developed independently in~\cite{ahmad181}.
\begin{algo}
\label{algo:main:1}
Approximate the expansion functions $d_k(\theta)$ in \eqref{eq:main:xiexpansion} in the points $\theta_{j,n_1}$ defined in \eqref{eq:introduction:thetagrid}.
\begin{enumerate}
\item Choose $n_1$ and $\alpha$.
\item For each $k=1,\ldots\alpha$
\begin{enumerate}
\item Compute eigenvalues $\lambda_j(A_{n_k})$ where, for example,
\begin{itemize}
\item $A_{n_k}=T_{n_k}(f)$,
\item $A_{n_k}=T_{n_k}(f)+R_{n_k}+N_{n_k}$,
\item $A_{n_k}=(T_{n_k}(b))^{-1}T_{n_k}(a)$, $f=a/b$,
\end{itemize}
and $n_k=2^{k-1}(n_1+1)-1$, $R_n$ is a low-rank matrix and $N_n$ is a small-norm matrix, $a,b$ are real-valued cosine polynomials, with $b$ non-negative and not identically zero; see~\cite{garoni171};
\item Define the indices $j_k=2^{k-1}j_1$, where $j_1=\{1,\ldots,n_1\}$;
\item Compute the grid that gives the exact eigenvalues by finding the roots $\xi_{j_k,n_k}$ such that 
\begin{align}
f(\xi_{j_k,n_k})-\lambda_{j_k}(A_{n_k})=0;\nonumber
\end{align}
\item Store the errors in a matrix $(\mathbf{E})_{k,j_k}=\xi_{j_k,n_k}-\theta_{j_1,n_1}$;
\end{enumerate}
\item Create the Vandermonde matrix $(\mathbf{V})_{i,j}=h_i^j$, where $h_i=1/(n_i+1)$, where $i,j=1,\ldots,\alpha$;
\item Compute the $\mathbf{D}$ matrix by solving $\mathbf{V}\mathbf{D}=\mathbf{E}$, that is, $\mathbf{D}=\mathbf{V}\backslash \mathbf{E}$.
\end{enumerate}
\end{algo}
\noindent Once we have our $\mathbf{D}$ matrix, where $\tilde{d}_k(\theta_{j,n_1})=(\mathbf{D})_{k,j_1}$, we can utilize them to reconstruct approximately the ``perfect'' grid $\xi_{j,n}$ for any $n$. Hence, for an arbitrary $n$, the full spectrum om $A_n$ is then approximated by \eqref{eq:main:lambdaapproxxiexpansion}, by first using a modified version of  the interpolation--extrapolation scheme of~\cite{ekstrom183}, to get $\tilde{d}_k(\theta_{j,n})$. The procedure is described in Algorithm~\ref{algo:main:2}. 
\begin{algo}
 \label{algo:main:2}
 Compute approximation $\tilde{d}_k(\theta_{j,n})$ from $\tilde{d}_k(\theta_{j,n_1})$, and approximate $\lambda_j(A_n)$.
\begin{enumerate}
\item Interpolate--extrapolate $\tilde{d}_k(\theta_{j,n_1})$ to $\tilde{d}_k(\theta_{j,n})$ as is done for $\tilde{c}_k(\theta_{j,n_1})$ in \cite{ekstrom183};
\item Compute $\tilde{\xi}_{j,n}^{(\alpha)}=\theta_{j,n}+\sum_{k=1}^\alpha\tilde{d}_k(\theta_{j,n})h^k$, where $h=1/(n+1)$;
\item Compute the eigenvalue approximations $\tilde{\lambda}_{j}^{(\xi,\alpha)}(A_{n})=f(\tilde{\xi}_{j,n}^{(\alpha)})$.
\end{enumerate}
\end{algo}

\section{Numerical Experiments}
\label{sec:numerical}
In this section we consider three different examples, numerically supporting Fact~\ref{fact:main:1}.
In Example~\ref{exmp:numerical:laplace} we study the Toeplitz matrix associated with the second order finite difference discretization of the Laplacian, with imposed Neumann--Dirichlet boundary conditions. 
We examine in Example~\ref{exmp:numerical:bilaplace}  the grid expansion \eqref{eq:main:xiapproxexpansion} for the Toeplitz matrix associated with second order finite difference discretization of the bi-Laplacian. 
Finally the grid expansion \eqref{eq:main:xiapproxexpansion} for a preconditioned matrix is presented in Example~\ref{exmp:numerical:precond}.
\begin{exmp}
\label{exmp:numerical:laplace}
We here study the symbol $f(\theta)=2-2\cos(\theta)$, and the perturbed generated matrix $A_n=T_n(f)+R_n$ of the form
\begin{align}
A_n=T_n(f)+R_n=
\left[
\begin{array}{rrrrrrrrr}
2&-1\\
-1&2&-1\\
&\ddots&\ddots&\ddots\\
&&\ddots&\ddots&-1\\
&&&-1&2
\end{array}
\right]
+
\left[
\begin{array}{rrrrrrrrr}
-1&\\
&\phantom{2}&\\
&\phantom{\ddots}&\phantom{\ddots}&\\
&&\phantom{\ddots}&&\\
&&&&
\end{array}
\right]
=
\left[
\begin{array}{rrrrrrrrr}
1&-1\\
-1&2&-1\\
&\ddots&\ddots&\ddots\\
&&\ddots&\ddots&-1\\
&&&-1&2
\end{array}
\right].\nonumber
\end{align}
This matrix corresponds to the normalized discretization of second order finite differences och the second derivative operator (Laplacian), with imposed mixed Neumann--Dirichlet boundary conditions.
It is known that the eigenvalues of $A_n$ are exactly expressed by
\begin{align}
\lambda_j(A_n)=f(\xi_{j,n}),\quad \xi_{j,n}=\frac{(j-1/2)\pi}{n+1/2},\quad j=1,\ldots,n.\nonumber
\end{align}
We thus have
\begin{align}
E_{j,n}^\xi&=\xi_{j,n}-\theta_{j,n}
=\frac{(j-1/2)\pi}{n+1/2}-\frac{j\pi}{n+1}
=\frac{1}{2n+1}\left(\theta_{j,n}-\pi\right),\nonumber
\end{align}
and then
\begin{align}
E_{j,n}^\xi=\sum_{k=1}^\infty d_k(\theta_{j,n})h^k,\quad h=1/(n+1),\nonumber
\end{align}
where
\begin{align}
d_k(\theta_{j,n})
=\frac{1}{2^k}\left(\theta_{j,n}-\pi\right).\nonumber
\end{align}
In Figure~\ref{fig:numerical:laplace} we present $\tilde{d}_k(\theta_{j,n_1})$, computed with Algorithm~\ref{algo:main:1}, and $d_k(\theta_{j,n_1})$ for $k=1,\ldots \alpha$, for $n_1=100, \alpha=4$. The overlaps of the numerical and theoretical function values are good.
\begin{figure}[!ht]
\centering
\includegraphics[width=0.48\textwidth]{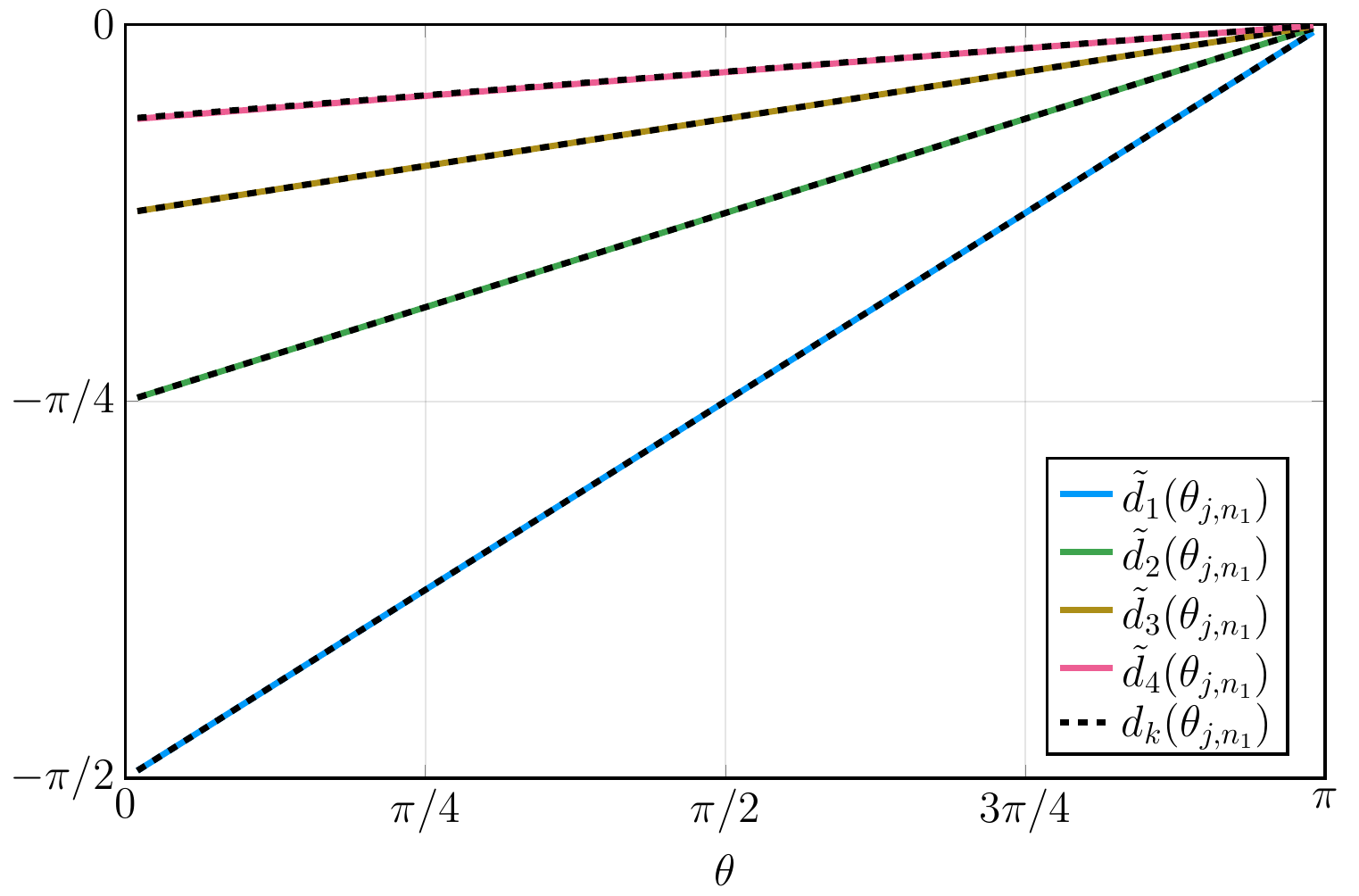}
\caption{Example \ref{exmp:numerical:laplace}: Expansion of the grid $\xi_{j,n_1}$ associated with the matrix $A_{n_1}=T_{n_1}(f)+R_{n_1}$, where $f(\theta)=2-2\cos(\theta)$ and $R_n$ is a small rank perturbation (coming from a Neumann boundary condition). The functions $d_k(\theta)$ are well approximated by $\tilde{d}_k(\theta)$.}
\label{fig:numerical:laplace}
\end{figure}
\end{exmp} 

\begin{exmp}
\label{exmp:numerical:bilaplace}
In this example we return to the symbol $f$ discussed in Examples~\ref{exmp:main:errorlambdavsxi} and \ref{exmp:main:errorxi}, namely,
\begin{align}
f(\theta)=(2-2\cos(\theta))^2=6-8\cos(\theta)+2\cos(2\theta),\nonumber
\end{align}
and employ Algorithms~\ref{algo:main:1} and \ref{algo:main:2} and study the results.

In Figure~\ref{fig:numerical:bilaplaceexpansion} we show the computed expansion functions $\tilde{d}_k(\theta_{j,n_1})$, for $n_1=100$ and $\alpha=2$ (left panel) and $\alpha=3$ (right panel). We note that for example $\tilde{d}_2(\theta_{1,n_1})$ and $\tilde{d}_2(\theta_{2,n_1})$ (red line, closest to $\theta=0$) is different in the two panels. Also $\tilde{d}_3(\theta_{j,n_1})$ behaves erratic a few approximations close to $\theta=0$.
This effect is due to the fact that the symbol $f$ violates the simple-loop conditions, discussed in detail in~\cite{barrera181} and also in~\cite{ekstrom184}. However, this has low impact in the current numerical setting, and also we employ a modification of the extrapolation in~\cite{ekstrom183} by simply ignoring user specified indices for the different $k=1,\ldots,\alpha$.

\begin{figure}[!ht]
\centering
\includegraphics[width=0.48\textwidth]{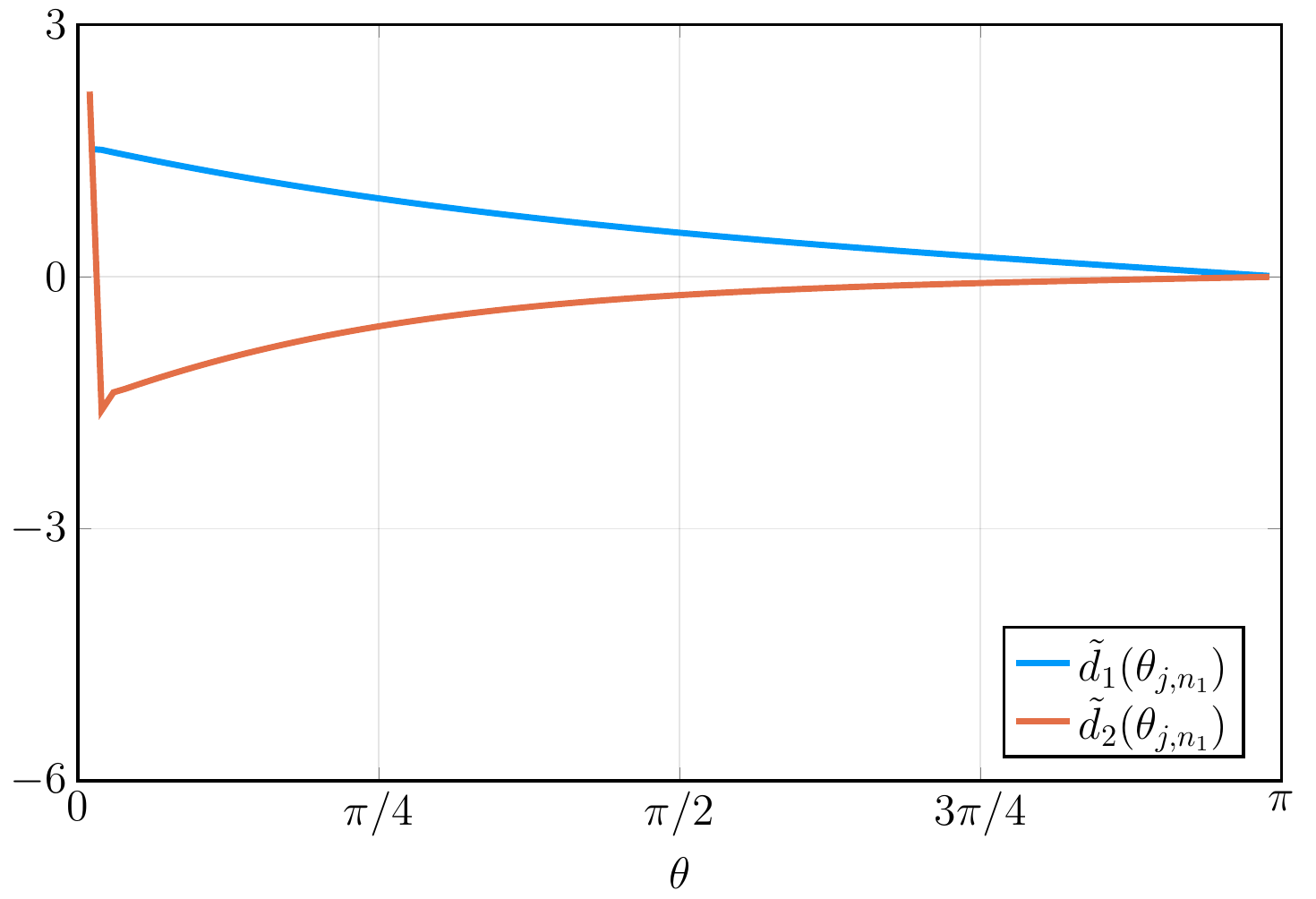}
\includegraphics[width=0.48\textwidth]{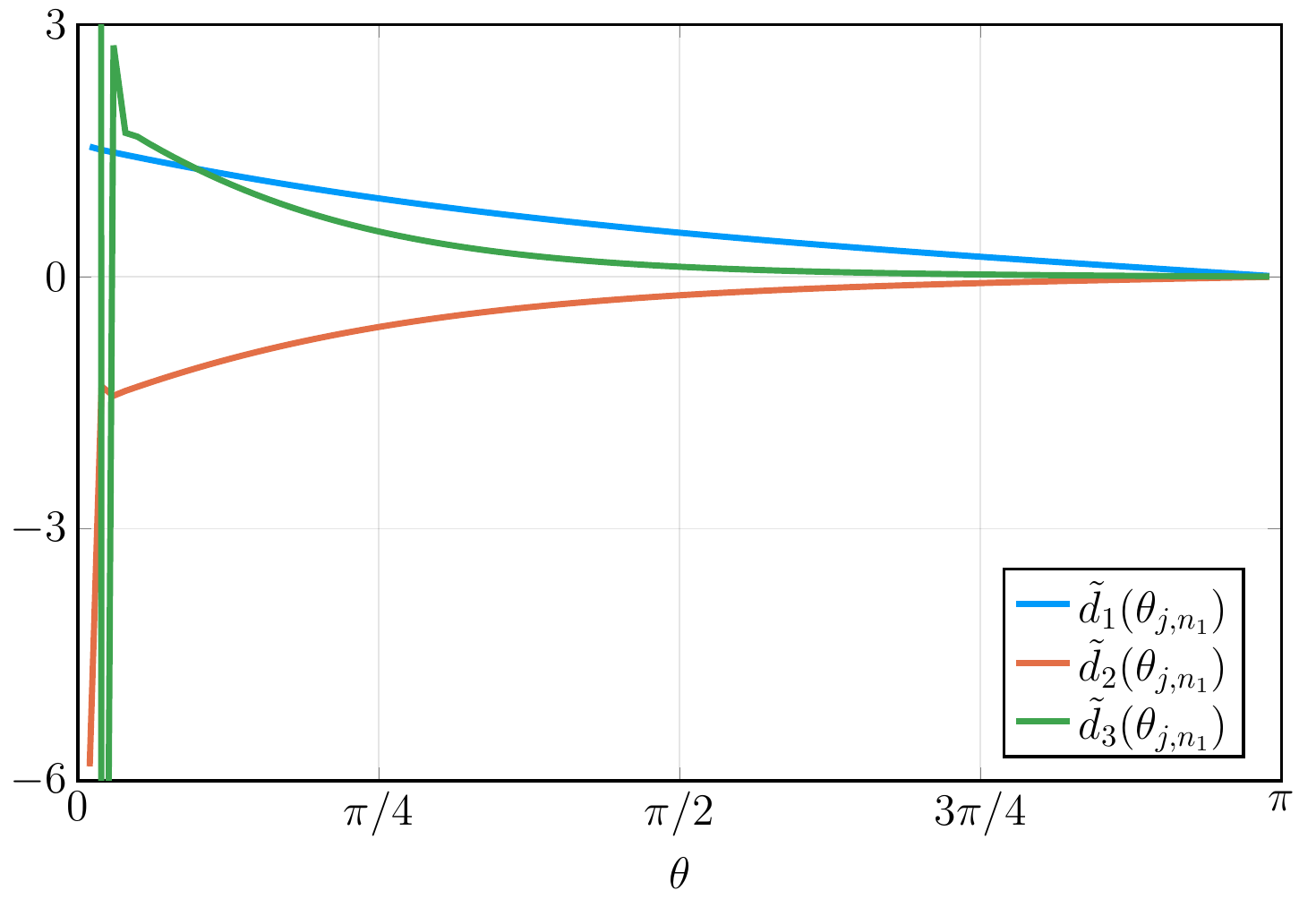}
\caption{Example~\ref{exmp:numerical:bilaplace}: Computed expansion function approximations $\tilde{d}_k(\theta_{j,n_1})$, $k=1,\ldots,\alpha$ for $\alpha=2$ (left) and $\alpha=3$ (right), and $n_1=100$. Note the erratic behavior for $\tilde{d}_2(\theta_{j,n_1})$ and $\tilde{d}_3(\theta_{j,n_1})$ close to $\theta=0$.}
\label{fig:numerical:bilaplaceexpansion}
\end{figure}
In the left panel of Figure~\ref{fig:numerical:bilaplacereductiongriderror} we show the interpolation--extrapolated $\tilde{d}_k(\theta_{j,n})$ (dashed black lines) for $n_1=100,\alpha=3$, and $n=100000$. The following samplings were ignored in these computations: for $\tilde{d}_2(\theta_{j,n_1})$ indices $j=1,2$ and for $\tilde{d}_3(\theta_{j,n_1})$ indices $j=1,2,3$. In the right panel of Figure~\ref{fig:numerical:bilaplacereductiongriderror} is presented the absolute error $\log_{10}|E_{j,n}^\xi|$ of \eqref{eq:main:lambdaxie0} and the subsequent errors $\log_{10}|\tilde{E}_{j,n,\beta}^\xi|$ of \eqref{eq:main:errorxiapproxxiexpansion} for $\beta=1,2,3$. The perfect grid $\xi_{j,n}$ is very well approximated, with a slight reduction of accuracy close to $\theta=0$.
The parameters used for the computations in the right panel of Figure~\ref{fig:numerical:bilaplacereductiongriderror} (and also both panels of Figure~\ref{fig:numerical:bilaplacereductionlambdaerror}) are $n_1=1000, \alpha=3$, and $n=100000$. 
As expected the correction $\tilde{d}_3(\theta_{j,n})h^3$ to $\tilde{E}_{j,n,2}^\xi$ (to get $\tilde{E}_{j,n,3}^\xi$) has only negligible effect since $h^3=\mathcal{O}(10^{-15})$. Only every $50^{\mathrm{th}}$ (every $100^{\mathrm{th}}$ for $\tilde{E}_{j,n,3}^\xi$) error is plotted for improved clarity in the visualizations.
\begin{figure}[!ht] 
\centering
\includegraphics[width=0.47\textwidth]{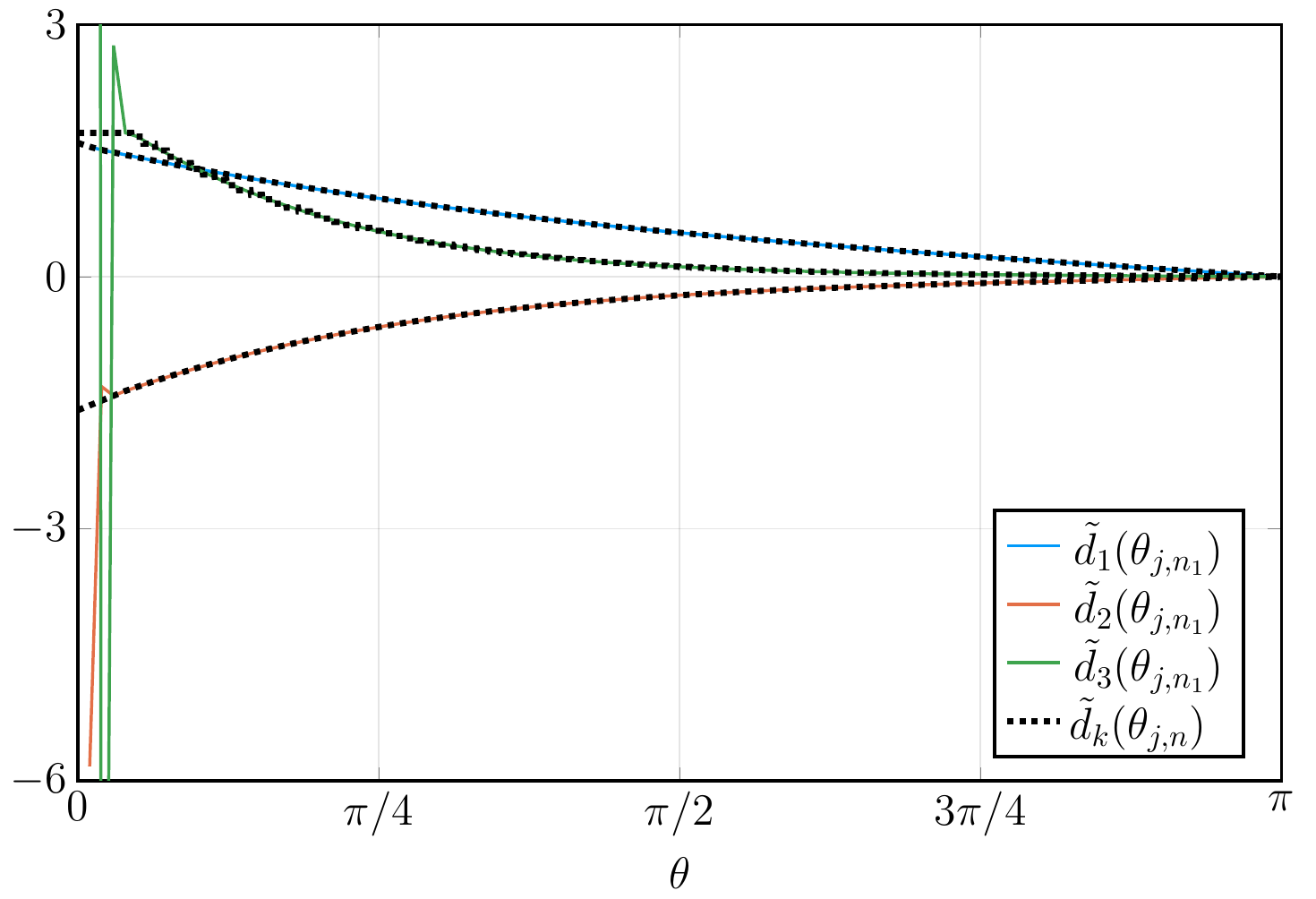}
\includegraphics[width=0.48\textwidth]{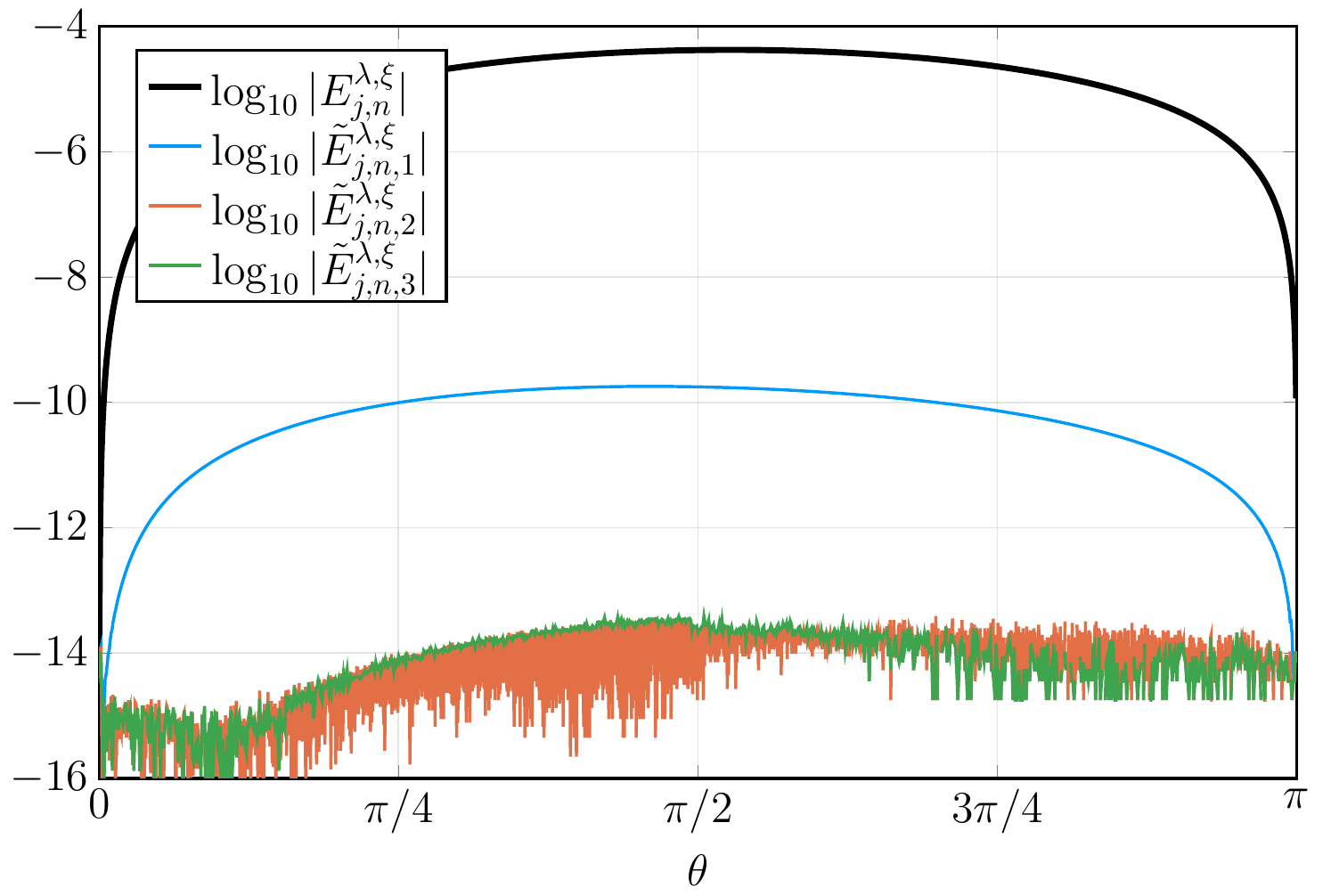}
\caption{Example~\ref{exmp:numerical:bilaplace}: Application of Algorithms~\ref{algo:main:1} and \ref{algo:main:2},. Left:  Approximations $\tilde{d}_k(\theta_{j,n_1})$ from Algorithm~\ref{algo:main:1} (colored lines) and $\tilde{d}_k(\theta_{j,n})$ from step 1 of of Algorithm~\ref{algo:main:2} (dashed black lines). Right: The error $E_{j,n}^\xi$ (and $\tilde{E}_{j,n,\beta}^\xi$, for $\beta=1,2,3)$.}
\label{fig:numerical:bilaplacereductiongriderror}
\end{figure}
In  Figure~\ref{fig:numerical:bilaplacereductionlambdaerror} we present the the errors $E_{j,n}^{\lambda,\theta}$ and $E_{j,n}^{\lambda,\xi}$ (and $\tilde{E}_{j,n,\beta}^{\lambda,\theta}$ and $\tilde{E}_{j,n,\beta}^{\lambda,\xi}$ for $\beta=1,2,3$) defined in \eqref{eq:introduction:lambdathetae0} and \eqref{eq:main:lambdaxie0} (and \eqref{eq:introduction:errorapproxexpansiontheta} and \eqref{eq:main:errorlambdaapproxxiexpansion}).
Thus, in the left panel of Figure~\ref{fig:numerical:bilaplacereductionlambdaerror} we show how well $\lambda_j(T_n(f))$ are approximated using the method in~\cite{ekstrom183} (with the addition of removal of erratic points in the interpolation--extrapolation stage). In the right panel of Figure~\ref{fig:numerical:bilaplacereductionlambdaerror} we present how well $\lambda_j(T_n(f))$ are approximated using Algorithms~\ref{algo:main:1} and \ref{algo:main:2} in this paper. We note that the result in the right panel is better overall for the whole spectrum, but especially close to $\theta=0$, and machine precision accuracy has almost been achieved. 
\begin{figure}[!ht] 
\centering
\includegraphics[width=0.48\textwidth]{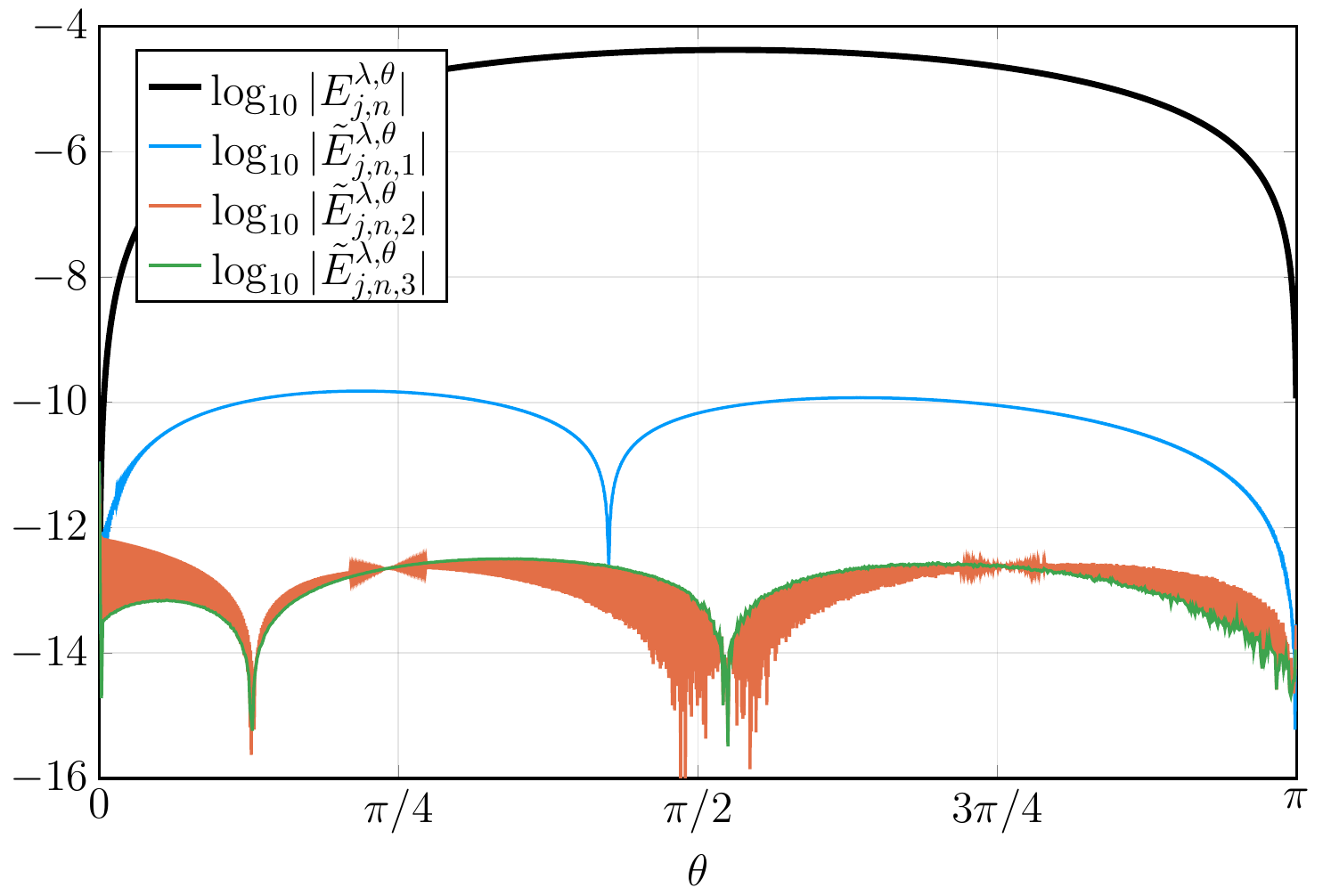}
\includegraphics[width=0.48\textwidth]{fig-numer-bilaplaceerrlx.pdf}
\caption{Example~\ref{exmp:numerical:bilaplace}: 
Errors of eigenvalue approximations. 
Left: Errors $\log_{10}|E_{j,n}^{\lambda,\theta}|$ (and $\log_{10}|\tilde{E}_{j,n,\beta}^{\lambda,\theta}|$, for $\beta=1,2,3$). Computations using~\cite{ekstrom183}.
Right: Errors $\log_{10}|E_{j,n}^{\lambda,\xi}|$ (and $\log_{10}|\tilde{E}_{j,n,\beta}^{\lambda,\xi}|$, for $\beta=1,2,3$). Computations using Algorithms~\ref{algo:main:1} and \ref{algo:main:2}. Note the overall lower error compared with the left panel, especially close to $\theta=0$.}
\label{fig:numerical:bilaplacereductionlambdaerror}
\end{figure}
\end{exmp}

\begin{exmp}
\label{exmp:numerical:precond}
Preconditioned matrices are important in many applications, and the study of their behavior is hence of importance. We here look at the matrix
\begin{align}
A_n&=T_n^{-1}(a)T_n(b),\nonumber\\
\{A_n\}&\sim_{\textsc{glt},\sigma,\lambda}f,\nonumber
\end{align}
where
\begin{align}
a(\theta)&=4-2\cos(\theta)-2\cos(2\theta)=(2-2\cos(\theta))(3+2\cos(\theta)),\nonumber\\
b(\theta)&=3+2\cos(\theta),\nonumber\\
f(\theta)&=a(\theta)/b(\theta)=2-2\cos(\theta),\nonumber
\end{align}
which is Example 1 from~\cite{ahmad171}. 

In the left panel of Figure~\ref{fig:numerical:precondsymbol} we show the symbols $a,b$, and $f=a/b$. In the right panel of Figure~\ref{fig:numerical:precondsymbol} is shown the absolute error when approximating $\xi_{j,n}$, that is, $\log_{10}|E_{j,n}^\xi|$ (and $\log_{10}|\tilde{E}_{j,n,\beta}^\xi|$ for $\beta=1,2,3$). Like in the right panel of Figure~\ref{fig:numerical:bilaplacereductiongriderror} we note a slight reduction of accuracy close to $\theta=0$, but here also close to $\theta=\pi$.

\begin{figure}[!ht] 
\centering
\includegraphics[width=0.47\textwidth]{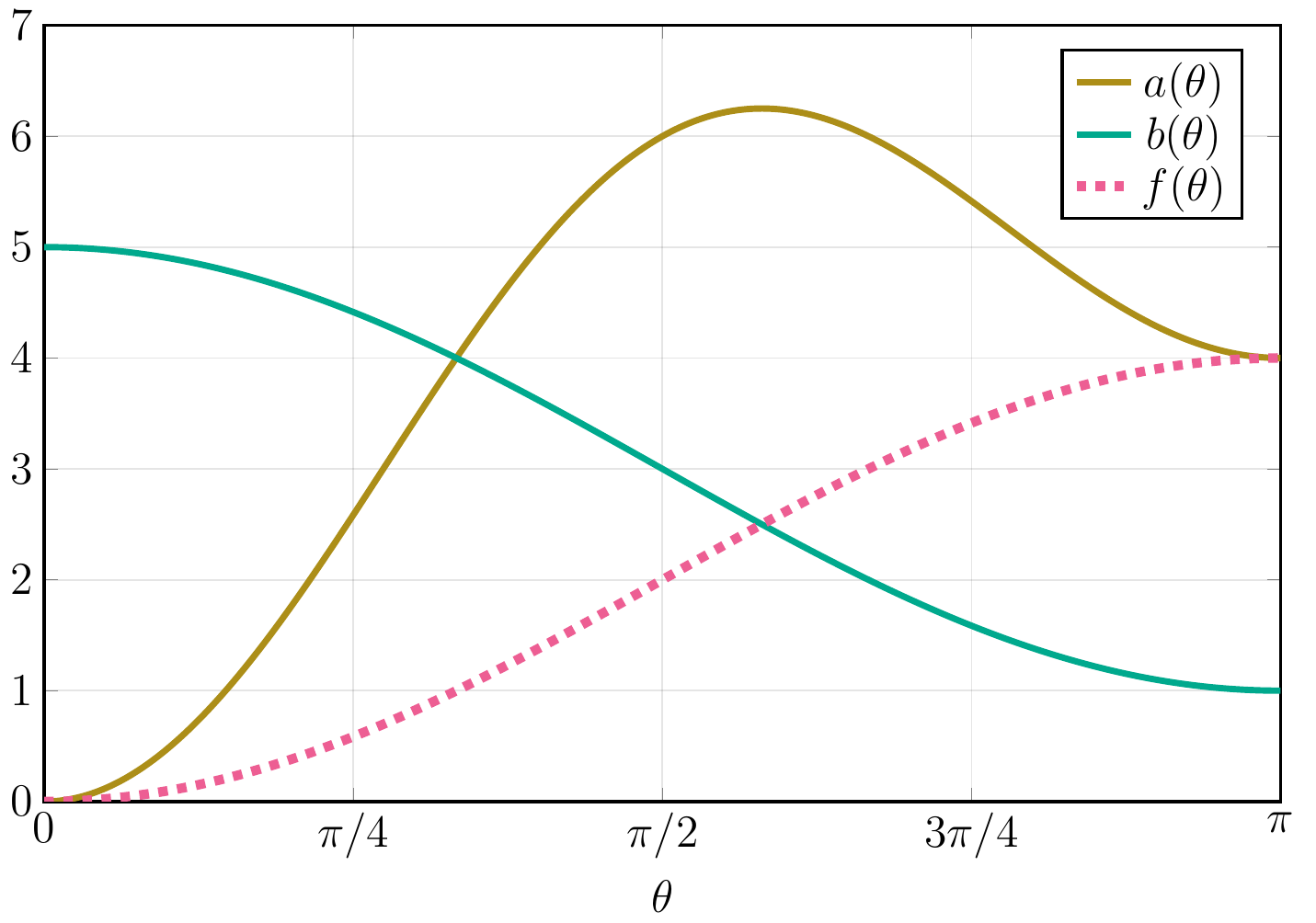}
\includegraphics[width=0.49\textwidth]{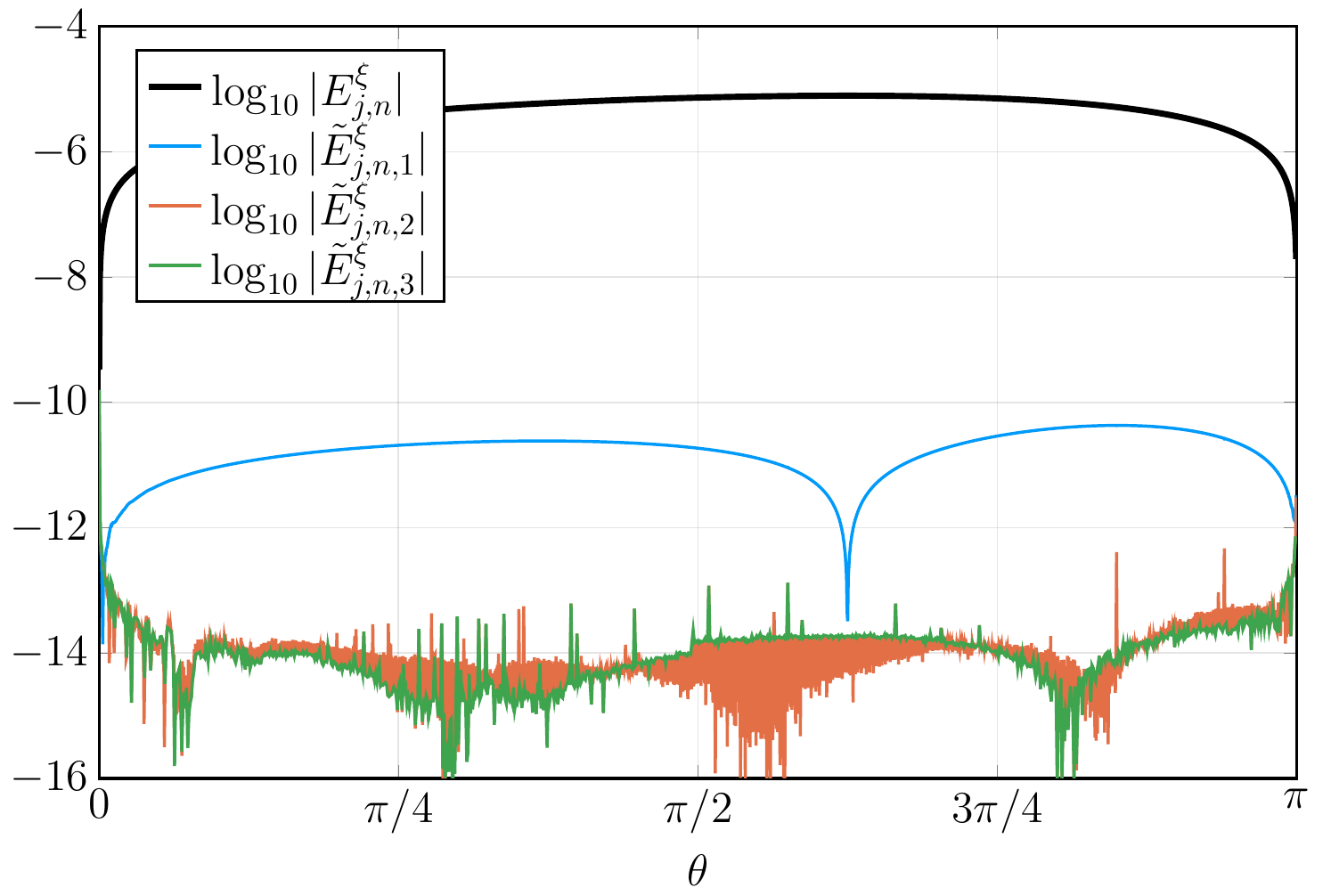}
\caption{Example~\ref{exmp:numerical:precond}: Symbols and errors approximating $\xi_{j,n}$. Left: Symbols $a,b$, and $f=a/b$. Right: $\log_{10}|E_{j,n}^\xi|$ (and $\log_{10}|\tilde{E}_{j,n,\beta}^\xi|$ for $\beta=1,2,3$ for $n_1=1000, \alpha=3$, and $n=100000$).}
\label{fig:numerical:precondsymbol}
\end{figure}

In Figure~\ref{fig:numerical:precondexpansion} we present both the computed expansion function $\tilde{c}_k(\theta_{j,n})$ (left) and $\tilde{d}_k(\theta_{j,n})$ (right). We note the similarity between the two expansions.
\begin{figure}[!ht] 
\centering
\includegraphics[width=0.47\textwidth]{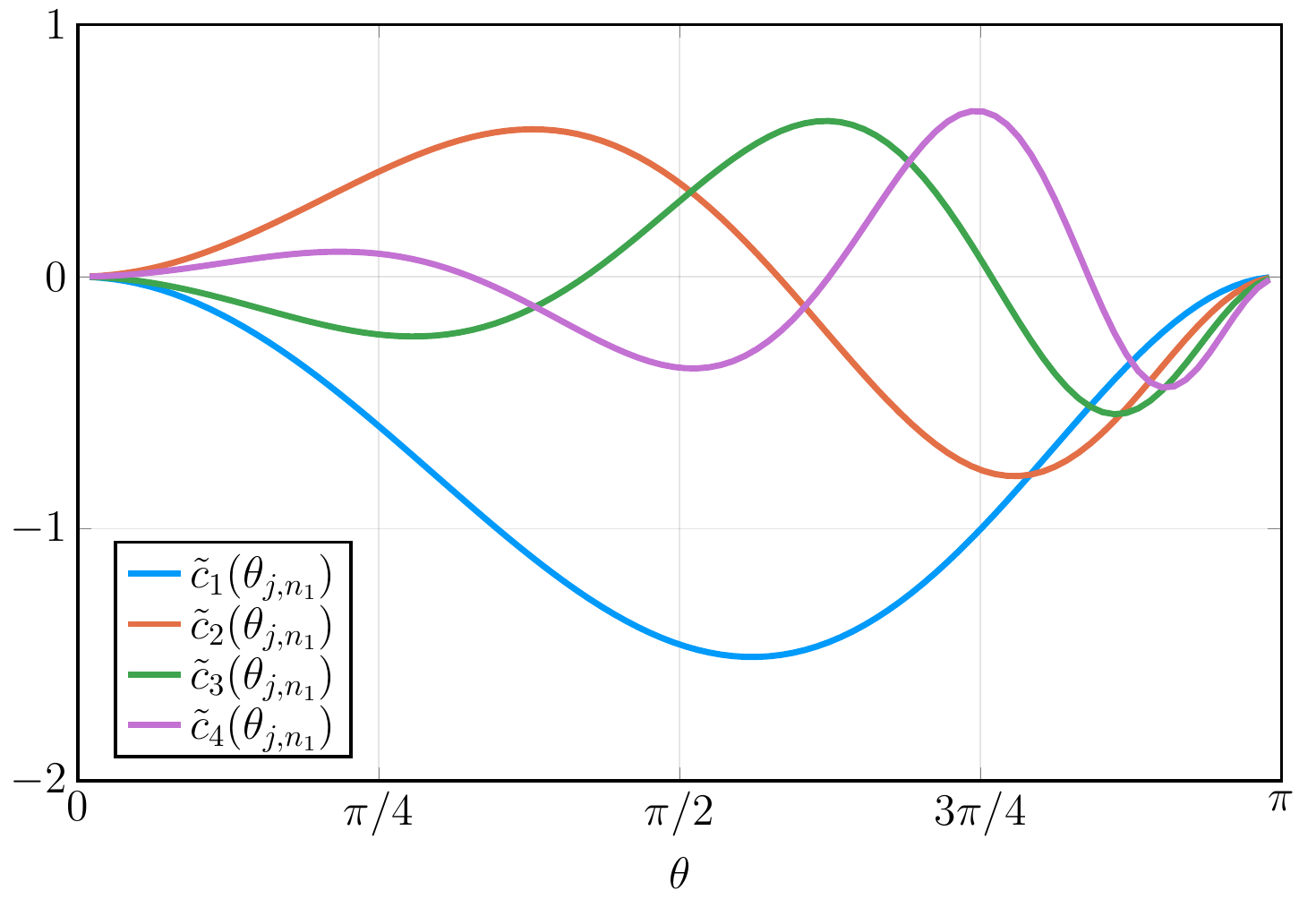}
\includegraphics[width=0.483\textwidth]{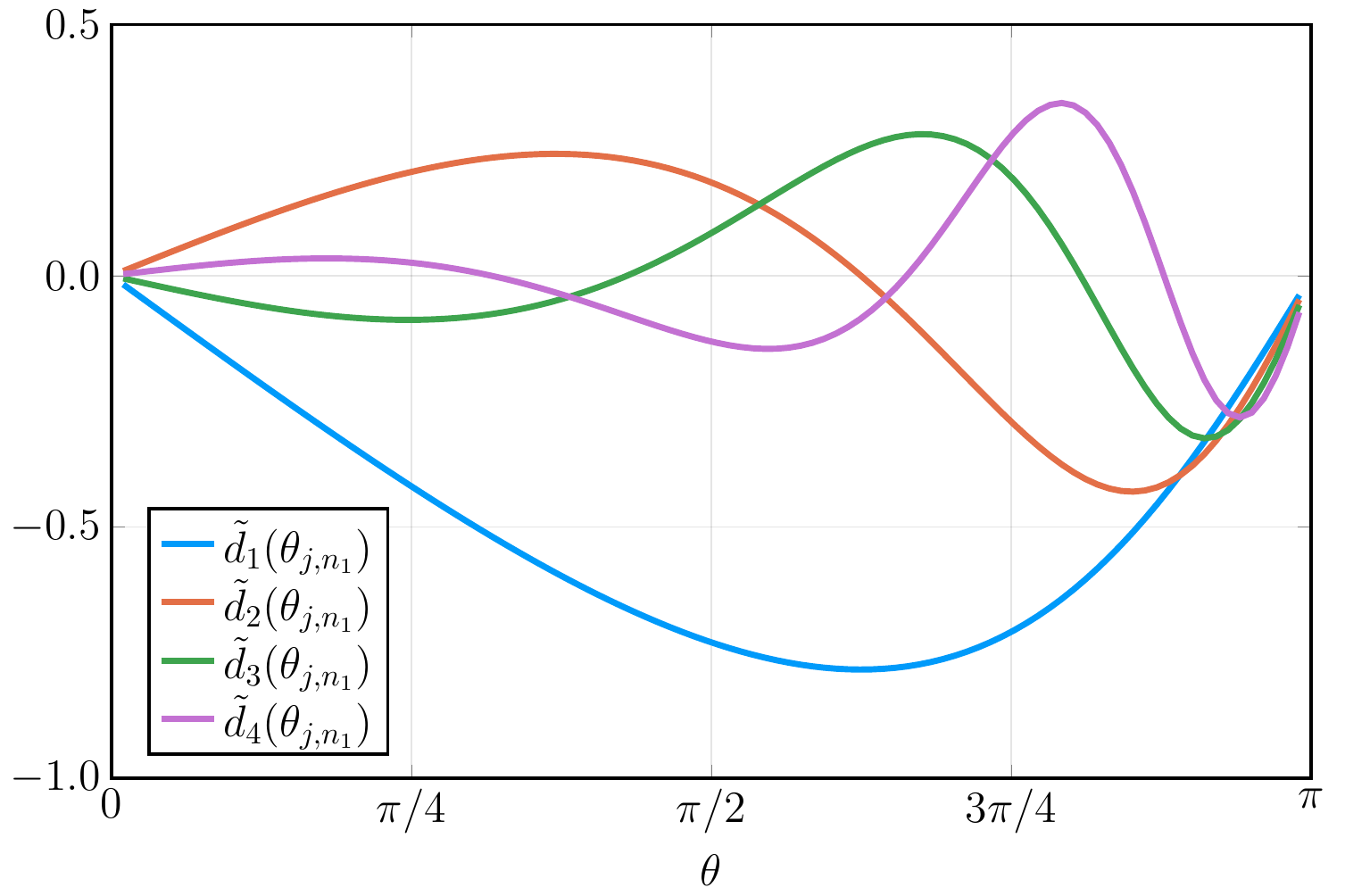}
\caption{Example~\ref{exmp:numerical:precond}: Expansions for preconditioned matrix, $n_1=100$ and $\alpha=4$. Left: Expansion $\tilde{c}_k(\theta_{j,n_1})$, for $k=1,\ldots,\alpha$. Right: Expansion $\tilde{d}_k(\theta_{j,n_1})$, for $k=1,\ldots,4$.}
\label{fig:numerical:precondexpansion}
\end{figure}
Finally in Figure~\ref{fig:numerical:precondlambdaerror} we present the absolute errors when approximating $\lambda_j(T_n^{-1}(b)T_n(a))$, that is, $\log_{10}|E_{j,n}^{\lambda,\theta}|$ (left panel) and $\log_{10}|E_{j,n}^{\lambda,\xi}|$ (right panel) (and $\log_{10}|\tilde{E}_{j,n,\beta}^{\lambda,\theta}|$ and $\log_{10}|\tilde{E}_{j,n,\beta}^{\lambda,\xi}|$ for $\beta=1,2,3$). 
As is noted in Figure~\ref{fig:numerical:bilaplacereductionlambdaerror}, the Algorithms~\ref{algo:main:1} and \ref{algo:main:2} generally perform better than the algorithm in~\cite{ekstrom183}.
The additional cost for Algorithms~\ref{algo:main:1} and \ref{algo:main:2}, compared with \cite{ekstrom183}, is negligible. It consists only of finding the $\alpha n_1$ roots to $f(\theta)-\lambda_{j_k}(A_{n_k})$, for $k=1,\ldots,\alpha$. 
\begin{figure}[!ht] 
\centering
\includegraphics[width=0.48\textwidth]{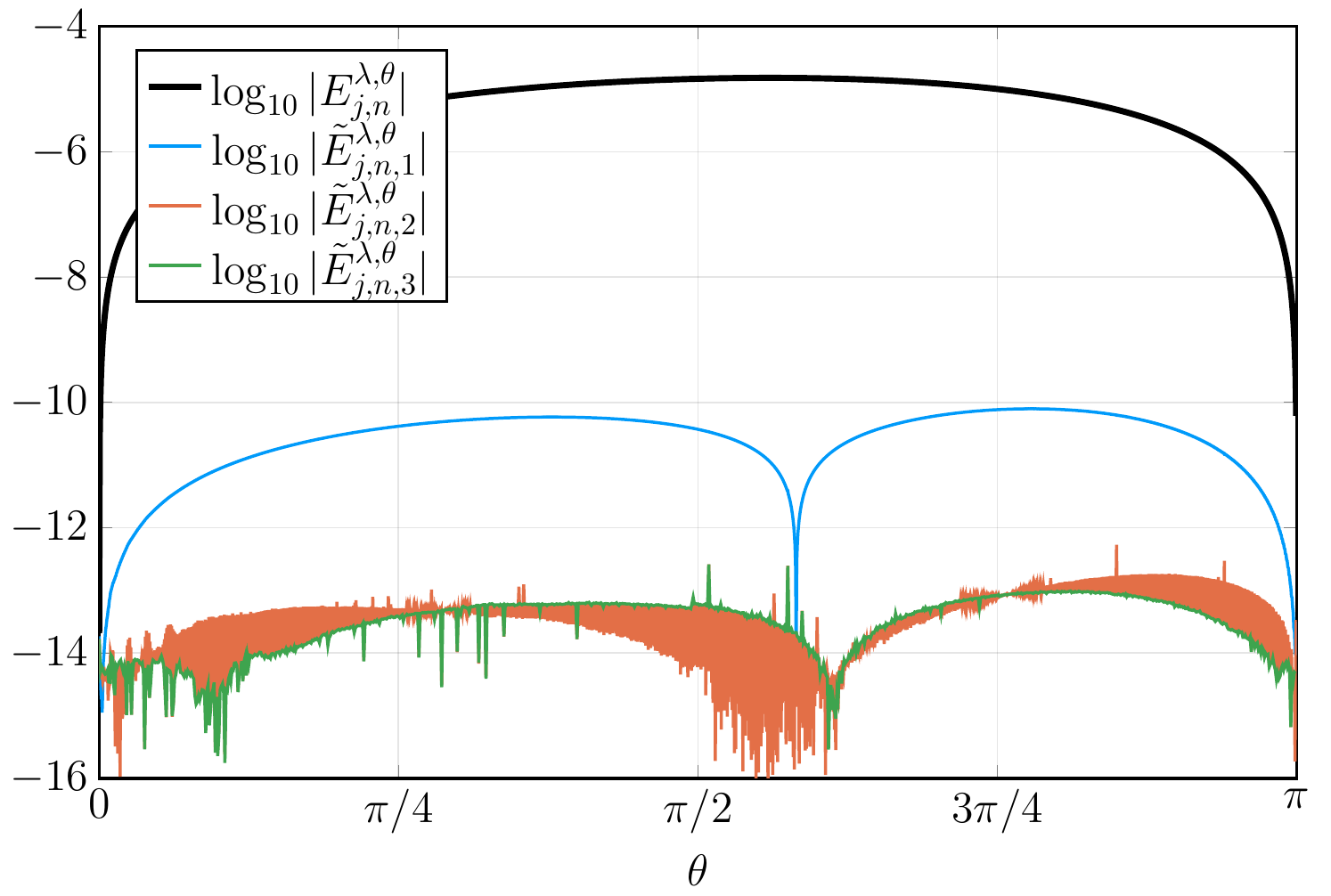}
\includegraphics[width=0.48\textwidth]{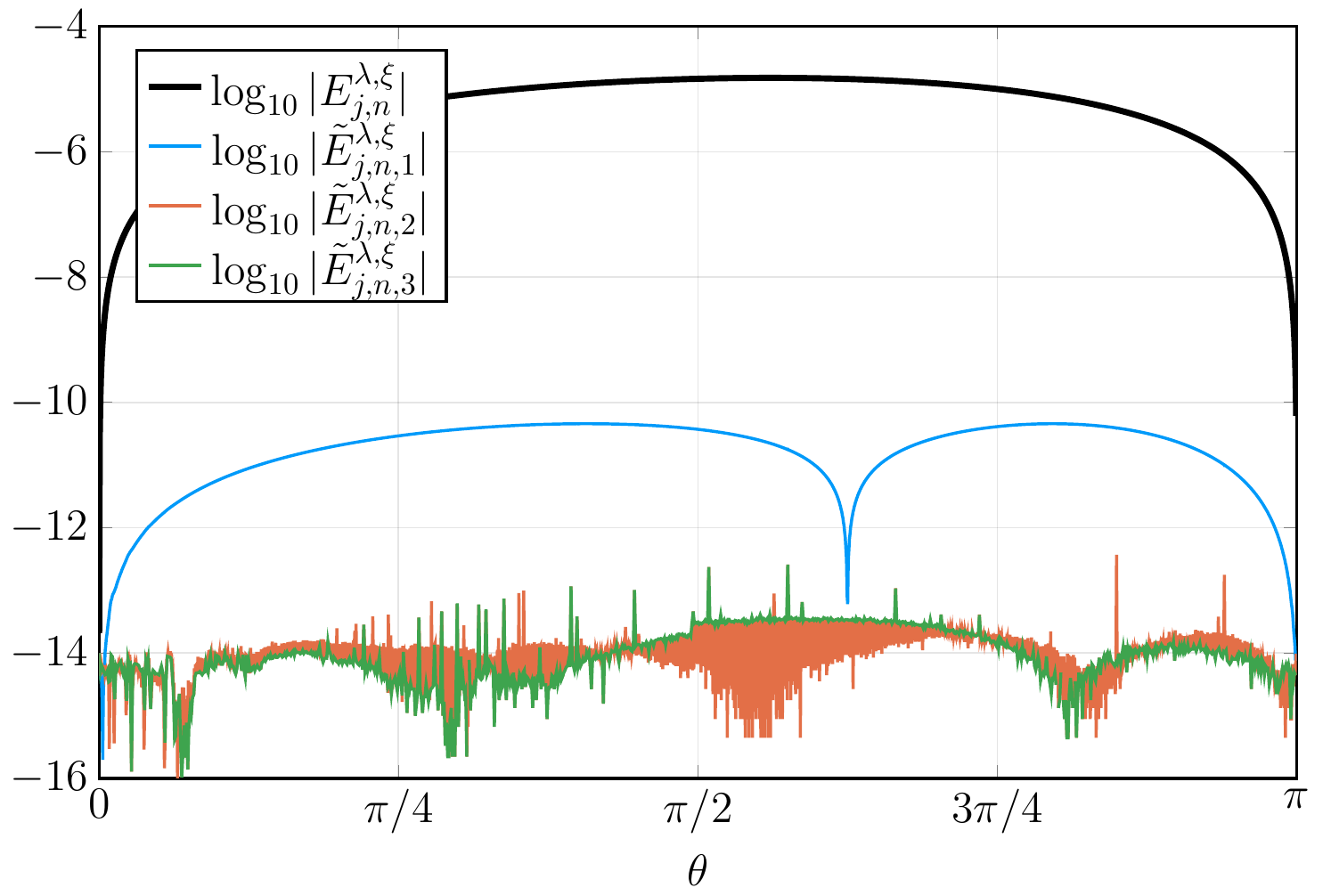}
\caption{Example~\ref{exmp:numerical:precond}: 
Errors of eigenvalue approximations. 
Left: Errors $\log_{10}|E_{j,n}^{\lambda,\theta}|$ (and $\log_{10}|\tilde{E}_{j,n,\beta}^{\lambda,\theta}|$, for $\beta=1,2,3$). Computations using~\cite{ekstrom183}.
Right: Errors $\log_{10}|E_{j,n}^{\lambda,\xi}|$ (and $\log_{10}|\tilde{E}_{j,n,\beta}^{\lambda,\xi}|$, for $\beta=1,2,3$). Computations using Algorithms~\ref{algo:main:1} and \ref{algo:main:2}. Note the overall lower error compared with the left panel.}
\label{fig:numerical:precondlambdaerror}
\end{figure}
\end{exmp}
 
\section{Conclusions}
In this paper we have studied the expansion of the grid error $E_{j,n}^\xi$ in $\xi_{j,n}=\theta_{j,n}+E_{j,n}^\xi$, where $\xi_{j,n}$ is the ``perfect'' grid associated with a sequence of matrices $\{A_n\}_n$, in the sense that $\lambda_j(A_n)=f(\xi_{j,n})$. Here $f$ is the symbol such that $\{A_n\}_n\sim_{\textsc{glt},\sigma,\lambda}f$. We present an Algorithm~\ref{algo:main:1} to approximate the expansion, proposed in Fact~\ref{fact:main:1}, of the error $E_{j,n}^\xi$ in \eqref{eq:main:xiexpansion}.
In Algorithm~\ref{algo:main:2} we show how to interpolate--extrapolate the data from Algorithm~\ref{algo:main:1}, and use this information to approximate the eigenvalues of $A_n$ for arbitrary $n$.
By the numerical experiments in Section~\ref{sec:numerical} we have shown evidence of Fact~\ref{fact:main:1}, in a number of practical examples, by also emphasizing the superiority of the new matrix-less method seen in Examples~\ref{exmp:numerical:bilaplace} and \ref{exmp:numerical:precond} when approximating the eigenvalues of a large matrix $A_n$.

The further study of the proposed expansion is warranted, since it might lead to new discoveries on the spectral behavior by applying perturbations of Toeplitz-like matrices, and also faster and more accurate methods to approximate the spectrum of Toeplitz-like matrices. For matrices generated by partially non-monotone symbols we expect similar behavior as for the expansion of the eigenvalue errors: in particular the conclusion is that we can not in general gain any information in the non-monotone part. The block and multilevel cases are also to be considered in future research.

\section{Acknowledgments}
The author would like to thank Stefano Serra-Capizzano for valuable discussions and suggestions during the preparation of this work. The author is financed by Athens University of Economics and Business.

\bibliography{References}

\begin{thebibliography}{10}

\bibitem{ahmad181}
{\sc F.~Ahmad}, {\em Equations and {S}ystems of {N}onlinear {E}quations: from
  high order {N}umerical Methods to fast {E}igensolvers for {S}tructured
  {M}atrices and {A}pplications}, PhD thesis, University of Insubria, 2018.

\bibitem{ahmad171}
{\sc F.~Ahmad, E.~S. Al-Aidarous, D.~A. Alrehaili, S.-E. Ekstr\"{o}m, I.~Furci,
  and S.~Serra-Capizzano}, {\em Are the eigenvalues of preconditioned banded
  symmetric {T}oeplitz matrices known in almost closed form?}, Numerical
  Algorithms, 78 (2017), pp.~867--893.

\bibitem{barrera181}
{\sc M.~Barrera, A.~B\"{o}ttcher, S.~M. Grudsky, and E.~A. Maximenko}, {\em
  Eigenvalues of even very nice {T}oeplitz matrices can be unexpectedly
  erratic}, in The Diversity and Beauty of Applied Operator Theory, Springer
  International Publishing, 2018, pp.~51--77.

\bibitem{bogoya151}
{\sc J.~Bogoya, A.~B\"{o}ttcher, S.~Grudsky, and E.~Maximenko}, {\em
  Eigenvalues of {H}ermitian {T}oeplitz matrices with smooth simple-loop
  symbols}, Journal of Mathematical Analysis and Applications, 422 (2015),
  pp.~1308--1334.

\bibitem{bogoya171}
{\sc J.~M. Bogoya, S.~M. Grudsky, and E.~A. Maximenko}, {\em Eigenvalues of
  {H}ermitian {T}oeplitz {M}atrices {G}enerated by {S}imple-loop {S}ymbols with
  {R}elaxed {S}moothness}, in Large Truncated Toeplitz Matrices, Toeplitz
  Operators, and Related Topics, Springer International Publishing, 2017,
  pp.~179--212.

\bibitem{bottcher101}
{\sc A.~B\"{o}ttcher, S.~Grudsky, and E.~Maksimenko}, {\em Inside the
  eigenvalues of certain {H}ermitian {T}oeplitz band matrices}, Journal of
  Computational and Applied Mathematics, 233 (2010), pp.~2245--2264.

\bibitem{bottcher991}
{\sc A.~B\"{o}ttcher and B.~Silbermann}, {\em Introduction to {L}arge
  {T}runcated {T}oeplitz {M}atrices}, Springer New York, 1999.

\bibitem{dibenedetto931}
{\sc F.~{Di Benedetto}, G.~Fiorentino, and S.~Serra}, {\em C. {G}.
  preconditioning for {T}oeplitz matrices}, Computers {\&} Mathematics with
  Applications, 25 (1993), pp.~35--45.

\bibitem{ekstrom185}
{\sc S.-E. Ekstr\"{o}m, I.~Furci, C.~Garoni, C.~Manni, S.~Serra-Capizzano, and
  H.~Speleers}, {\em Are the eigenvalues of the {B}-spline isogeometric
  analysis approximation of {$-\Delta u=\lambda u$} known in almost closed
  form?}, Numerical Linear Algebra with Applications, 25 (2018), p.~e2198.

\bibitem{ekstrom181}
{\sc S.-E. Ekstr\"{o}m, I.~Furci, and S.~Serra-Capizzano}, {\em Exact formulae
  and matrix-less eigensolvers for block banded symmetric {T}oeplitz matrices},
  {BIT} Numerical Mathematics, 58 (2018), pp.~937--968.

\bibitem{ekstrom183}
{\sc S.-E. Ekstr\"{o}m and C.~Garoni}, {\em A matrix-less and parallel
  interpolation--extrapolation algorithm for computing the eigenvalues of
  preconditioned banded symmetric {T}oeplitz matrices}, Numerical Algorithms,
  (2018).
\newblock https://doi.org/10.1007/s11075-018-0508-0, (in press).

\bibitem{ekstrom171}
{\sc S.-E. Ekstr\"{o}m, C.~Garoni, and S.~Serra-Capizzano}, {\em Are the
  {E}igenvalues of {B}anded {S}ymmetric {T}oeplitz {M}atrices {K}nown in
  {A}lmost {C}losed {F}orm?}, Experimental Mathematics, 27 (2018),
  pp.~478--487.

\bibitem{ekstrom184}
{\sc S.-E. Ekström}, {\em Matrix-{L}ess {M}ethods for {C}omputing
  {E}igenvalues of {L}arge {S}tructured {M}atrices}, PhD thesis, Uppsala
  University, Uppsala: Acta Universitatis Upsaliensis, 2018.

\bibitem{garoni171}
{\sc C.~Garoni and S.~Serra-Capizzano}, {\em Generalized {L}ocally {T}oeplitz
  {S}equences: {T}heory and {A}pplications ({V}olume 1)}, Springer
  International Publishing, 2017.

\end{thebibliography}
\bibliographystyle{siam}
\end{document}